\documentclass[11pt]{amsart}
\usepackage{graphicx, color}
\usepackage{amscd}
\usepackage{amsmath}
\usepackage{amsfonts}
\usepackage{amssymb}
\usepackage{mathrsfs}

\newtheorem{theorem}{Theorem}

\newtheorem{remark}[theorem]{Remark}

\numberwithin{equation}{section}

\renewcommand{\leq}{\leqslant}

\renewcommand{\geq}{\geqslant}
\begin{document}
\title[Hessian equations and systems with weights]{A necessary and a
sufficient condition for the existence of the positive radial solutions to
Hessian equations and systems with weights}
\author[D.-P. Covei]{Dragos-Patru Covei}
\address{Department of Applied Mathematics, The Bucharest University of
Economic Studies, Piata Romana, 1st district, postal code: 010374, postal
office: 22, Romania}
\email{\texttt{coveidragos@yahoo.com}}

\begin{abstract}
{\footnotesize In this article we consider the existence} {\footnotesize of
positive radial solutions for\textbf{\ }Hessian equations and systems with
weights and we give a necessary condition as well as a sufficient condition
for a positive radial solution to be large. The method of proving theorems
is essentially based on a successive approximation. Our results complete and
improve a recently work published by Zhang and Zhou (\textit{Existence of
entire positive k-convex radial solutions to Hessian equations and systems
with weights}, Applied Mathematics Letters, Volume 50, December 2015, Pages
48--55).}
\end{abstract}

\maketitle
\tableofcontents




\section{Introduction}

Let $D^{2}u$ be the Hessian matrix of a $C^{2}$ (i.e., a twice continuously
differentiable) function $u$ defined over $\mathbb{R}^{N}$ ($N\geq 3$) and $%
\lambda \left( D^{2}u\right) =\left( \lambda _{1},...,\lambda _{N}\right) $
the vector of eigenvalues of $D^{2}u$. For $k=1,2,...,N$ is defined the $k$%
-Hessian operator as follows%
\begin{equation*}
S_{k}\left( \lambda \left( D^{2}u\right) \right) =\underset{1\leq
i_{1}<...<i_{k}\leq N}{\sum }\lambda _{i_{1}}\cdot ...\cdot \lambda _{i_{k}}
\end{equation*}%
i.e., it is the $k^{th}$ elementary symmetric polynomial of the Hessian
matrix of $u$. In other words, $S_{k}\left( \lambda \left( D^{2}u\right)
\right) $ it is the sum of all $k\times k$ principal minors of the Hessian
matrix $D^{2}u$ and so is a second order differential operator, which may
also be called the $k$-trace of $D^{2}u$. Especially, it is easily to see
that the $N$-Hessian is the Monge-Amp\'{e}re operator and that the $1$%
-Hessian is the well known classical Laplace operator. Hence, the $k$%
-Hessian operators form a discrete collection of partial differential
operators which includes both the Laplace and the Monge-Amp\'{e}re operator.

In this paper we study the existence of radial solutions for the following
Hessian equation 
\begin{equation}
S_{k}^{1/k}\left( \lambda \left( D^{2}u\right) \right) =p\left( \left\vert
x\right\vert \right) h\left( u\right) \text{ in }\mathbb{R}^{N}\text{ ,}
\label{11}
\end{equation}%
and system%
\begin{equation}
\left\{ 
\begin{array}{c}
S_{k}^{1/k}\left( \lambda \left( D^{2}u\right) \right) =p\left( \left\vert
x\right\vert \right) f\left( u,v\right) \text{ in }\mathbb{R}^{N}\text{ ,}
\\ 
S_{k}^{1/k}\left( \lambda \left( D^{2}v\right) \right) =q\left( \left\vert
x\right\vert \right) g\left( u,v\right) \text{ in }\mathbb{R}^{N}\text{ ,}%
\end{array}%
\right.  \label{11s}
\end{equation}%
where $k\in \left\{ 1,2,...,N\right\} $, the continuous functions $p$, $q$%
\textit{\ }$:\left[ 0,\infty \right) \rightarrow \left( 0,\infty \right) $, $%
h:\left[ 0,\infty \right) \rightarrow \left[ 0,\infty \right) $ and $f$, $g:%
\left[ 0,\infty \right) \times \left[ 0,\infty \right) \rightarrow \left[
0,\infty \right) $ satisfy some of the conditions:

(P1)\quad $p$, $q$ is a spherically symmetric function (i.e.\textit{\ }$%
p\left( x\right) =p\left( \left\vert x\right\vert \right) $, $q\left(
x\right) =q\left( \left\vert x\right\vert \right) );$

(P2)\quad $r^{N+\frac{N}{k}-2}p^{k}\left( r\right) $\textit{\ }is
nondecreasing for large $r$;

(P3)\quad $r^{N+\frac{N}{k}-2}\left[ p^{k}\left( r\right) +q^{k}\left(
r\right) \right] $\textit{\ }is nondecreasing for large $r$\textit{;}

(C1)\quad $h$ is monotone non-decreasing, $h(0)=0$ and $h\left( s\right) >0$
for all $s>0$;

(C2)\quad $f$, $g$ are monotone non-decreasing in each variable, $%
f(0,0)=g\left( 0,0\right) =0$ and $f(s,t)>0,$ $g\left( s,t\right) >0$ for
all $s,t>0$;

(C3)\quad $\int_{1}^{\infty }\frac{1}{\sqrt[k+1]{\left( k+1\right) H(t)}}%
dt=\infty $ for\ $H(t)=\int_{0}^{t}h^{k}(z)dz;$

(C4)\quad $\int_{1}^{\infty }\frac{1}{\sqrt[k+1]{\left( k+1\right) F(t)}}%
dt=\infty $ for\ $F(t)=\int_{0}^{t}\left( f^{k}(z,z)+g^{k}(z,z)\right) dz.$

The properties of the $k$-Hessian operator was well discussed in a numerous
papers written as a first author by Ivochkina (see \cite{III}-\cite{IF} and
others). Moreover, this operator appear as an object of investigation by
many remarkable geometers. For example, Viaclovsky (see \cite{V}, \cite{VI})
observed that the $k-$Hessian operator is an important class of fully
nonlinear operators which is closely related to a geometric problem of the
type (\ref{11}), where we cite the work of Bao-Ji-Li \cite{BAOII} for a more
detailed discussion. Moreover, equation (\ref{11}) arises via the study of
the quasilinear parabolic problem ( see for example the introduction of
Moll-Petitta \cite{MP}). In the present work we will limit ourselves to the
development of mathematical theory for (\ref{11}) and (\ref{11s}). The main
difficulty in investigating problems, such as (\ref{11}) or (\ref{11s}), in
which appear the $k$-Hessian operator is related to the fact that their
properties change depending on the subset of $C^{2}$ from where the solution
is taken. Our main objective here is to find functions in $C^{2}$ that are
strictly $k$-convex and verifies the problems (\ref{11}), (\ref{11s}), where
by strictly $k$-convex function $u$ we mean that all eigenvalues $\lambda
_{1},...,\lambda _{N}$ of the symmetric matrix $D^{2}u$ are in the so called
G\aa rdding open cone $\Gamma _{k}$ which is defined by%
\begin{equation*}
\Gamma _{k}\left( N\right) =\left\{ \lambda \in \mathbb{R}^{N}\left\vert
S_{1}\left( \lambda \right) >0,....,S_{k}\left( \lambda \right) >0\right.
\right\}
\end{equation*}%
In the next we adopt the notation from Bao-Li \cite{BAO} for the space of
all admissible functions 
\begin{equation*}
\Phi ^{k}\left( \mathbb{R}^{N}\right) :=\left\{ u\in C^{2}\left( \mathbb{R}%
^{N}\right) \left\vert \lambda \in \Gamma _{k}\left( N\right) \text{ for all 
}x\in \mathbb{R}^{N}\right. \right\} .
\end{equation*}%
In our direction, there are some recently papers resolving existence for
blow-up solutions of (\ref{11}) and (\ref{11s}). Here we wish to mention the
works of Bao-Ji-Li \cite{BAOII}, Jacobsen \cite{J}, Bao-Li \cite{BAO}, Lazer
and McKenna \cite[(the case $k=N$)]{LM}, Salani \cite{S} and Zhang-Zhou \cite%
{ZZ} which will be useful in our proofs. It is interesting to note that in
our results the dimension of the space $\mathbb{R}^{N}$ affect the
properties of the solution of the equation and system which in the case of
the classical Laplace operator and the Monge-Amp\'{e}re operator this
condition doesn't appear in any works.

Motivated by the recent work of Zhang-Zhou \cite{ZZ} we are interested in
proving the following theorems:

\begin{theorem}
\label{th1}Let $k\in \left\{ 1,2,...,\left[ N/2\right] \right\} $ if $N$ is
odd or $k\in \left\{ 1,2,...,\left[ N/2\right] -1\right\} $ if $N$ is even.
Suppose that \textrm{(P1), (P2)}, \textrm{(C1)}, \textrm{(C3)} are
satisfied. If there exists a positive number $\varepsilon $ such that 
\begin{equation}
\quad \text{ }\int_{0}^{\infty }t^{1+\varepsilon +\frac{2\left( k-1\right) }{%
k+1}}\left( p\left( t\right) \right) ^{\frac{2k}{k+1}}dt<\infty ,  \label{5}
\end{equation}%
then system \textrm{(\ref{11})} has a nonnegative nontrivial radial bounded
solution $u\in \Phi ^{k}\left( \mathbb{R}^{N}\right) $.
\end{theorem}

\begin{theorem}
\label{th2}\textit{If }$p$ satisfy \textit{\textrm{(P1)}} and\ $f$ \textit{%
satisfy \textrm{(C1)}, \textrm{(C3)}, then the problem \textrm{(\ref{11})}
has a nonnegative nontrivial entire radial solution }$u\in \Phi ^{k}\left( 
\mathbb{R}^{N}\right) $\textit{. Suppose furthermore that \textrm{(P2)}}
holds.\textit{\ }If $p$ satisfies 
\begin{equation}
\int_{0}^{\infty }\left( \frac{k}{t^{N-k}}\int_{0}^{t}\frac{s^{N-1}}{%
C_{N-1}^{k-1}}p^{k}\left( s\right) ds\right) ^{1/k}dt=\infty ,\text{ }
\label{12}
\end{equation}%
then any nonnegative nontrivial radial solution $u\in \Phi ^{k}\left( 
\mathbb{R}^{N}\right) $ of \textrm{(\ref{11})} is large. Conversely, \textit{%
if }\textrm{(\ref{11})} has a nonnegative entire large radial solution $u\in
\Phi ^{k}\left( \mathbb{R}^{N}\right) $, then one or both of the following 
\begin{equation}
\begin{array}{ll}
1. & \int_{0}^{\infty }t^{1+\varepsilon +\frac{2\left( k-1\right) }{k+1}%
}\left( p\left( t\right) \right) ^{\frac{2k}{k+1}}dr=\infty \text{ for every 
}\varepsilon >0\,; \\ 
2. & k\in \left\{ \left[ N/2\right] +1,...,N\right\} \text{ if }N\text{ is
odd or }k\in \left\{ \left[ N/2\right] ,...,N\right\} \text{ if }N\text{ is
even,}%
\end{array}
\label{13}
\end{equation}%
hold.
\end{theorem}

Regarding existence of solution to (\ref{11s}), we have the following
results.

\begin{theorem}
\label{th3}Let $k\in \left\{ 1,2,...,\left[ N/2\right] \right\} $ if $N$ is
odd or $k\in \left\{ 1,2,...,\left[ N/2\right] -1\right\} $ if $N$ is even.
Suppose that \textrm{(P1), (P3)}, \textrm{(C2)}, \textrm{(C4)} are
satisfied. If there exists a positive number $\varepsilon $ such that 
\begin{equation}
\text{ }\int_{0}^{\infty }t^{1+\varepsilon +\frac{2\left( k-1\right) }{k+1}%
}\left( p^{k}\left( t\right) +q^{k}\left( t\right) \right) ^{\frac{2}{k+1}%
}dt<\infty \text{ },  \label{5s}
\end{equation}%
then system \textrm{(\ref{11s})} has a nonnegative nontrivial radial bounded
solution $\left( u,v\right) \in \Phi ^{k}\left( \mathbb{R}^{N}\right) \times
\Phi ^{k}\left( \mathbb{R}^{N}\right) $.
\end{theorem}

\begin{theorem}
\label{th4}\textit{If }$p,$ $q$ satisfy \textit{\textrm{(P1)}} and\textit{\ }%
\ $f,$ $g$ \textit{satisfy \textrm{(C2)}, \textrm{(C4)}, then the problem 
\textrm{(\ref{11})} has a nonnegative nontrivial entire radial solution.
Suppose furthermore that \textrm{(P3)} holds. }If $p$ satisfies 
\begin{equation}
\int_{0}^{\infty }\left( \frac{k}{t^{N-k}}\int_{0}^{t}\frac{s^{N-1}}{%
C_{N-1}^{k-1}}p^{k}\left( s\right) ds\right) ^{1/k}dt=\infty \text{ and }%
\int_{0}^{\infty }\left( \frac{k}{t^{N-k}}\int_{0}^{t}\frac{s^{N-1}}{%
C_{N-1}^{k-1}}q^{k}\left( s\right) ds\right) ^{1/k}dt=\infty ,  \label{12s}
\end{equation}%
then any nonnegative nontrivial solution $\left( u,v\right) \in \Phi
^{k}\left( \mathbb{R}^{N}\right) \times \Phi ^{k}\left( \mathbb{R}%
^{N}\right) $ of \textrm{(\ref{11s})} is large. Conversely, if \textrm{(\ref%
{11s})} has a nonnegative entire large radial solution $\left( u,v\right)
\in \Phi ^{k}\left( \mathbb{R}^{N}\right) \times \Phi ^{k}\left( \mathbb{R}%
^{N}\right) $, then one or both of the following 
\begin{equation}
\begin{array}{ll}
1. & \int_{0}^{\infty }t^{1+\varepsilon +\frac{2\left( k-1\right) }{k+1}%
}\left( p^{k}\left( t\right) +q^{k}\left( t\right) \right) ^{\frac{2}{k+1}%
}dr=\infty \text{ for every }\varepsilon >0\,; \\ 
2. & k\in \left\{ \left[ N/2\right] +1,...,N\right\} \text{ if }N\text{ is
odd or }k\in \left\{ \left[ N/2\right] ,...,N\right\} \text{ if }N\text{ is
even,}%
\end{array}
\label{13s}
\end{equation}%
hold.
\end{theorem}

For the readers' convenience, we recall the radial form of the $k$-Hessian
operator.

\begin{remark}
(see, for example, \cite{BAO}, \cite{S}) If $u:\mathbb{R}^{N}\rightarrow 
\mathbb{R}$ is radially symmetric then a calculation show%
\begin{equation*}
S_{k}\left( \lambda \left( D^{2}u\left( r\right) \right) \right)
=r^{1-N}C_{N-1}^{k-1}\left[ \frac{r^{N-k}}{k}\left( u^{^{\prime }}\left(
r\right) \right) ^{k}\right] ^{\prime }\text{ ,}
\end{equation*}%
where the prime denotes differentiation with respect to $r=\left\vert
x\right\vert $ and $C_{N-1}^{k-1}=(N-1)!/\left[ (k-1)!(N-k)!\right] $.
\end{remark}

\section{Proofs of the main results}

In this section we give the proofs of Theorems \ref{th1} - \ref{th4}. The
main references for proving Theorems \ref{th1} - \ref{th2} is the work of
Lair \cite{LA} and Delano\"{e} \cite{DEL} see also Afrouzi-Shokooh \cite{AF}.

\subparagraph{Proof of the Theorem \protect\ref{th1}}

Assume that (\ref{5}) holds. We prove the existence of $w\in \Phi ^{k}\left( 
\mathbb{R}^{N}\right) $ to the problem

\begin{equation}
\begin{array}{l}
S_{k}^{1/k}\left( \lambda \left( D^{2}w\left( \left\vert x\right\vert
\right) \right) \right) =p\left( \left\vert x\right\vert \right) h\left(
w\left( \left\vert x\right\vert \right) \right) \text{ in }\mathbb{R}^{N}%
\text{. }%
\end{array}
\label{6}
\end{equation}%
Observe that we can rewrite (\ref{6}) as follows: 
\begin{equation*}
\left[ \frac{r^{N-k}}{k}\left( w^{^{\prime }}(r)\right) ^{k}\right] ^{\prime
}=\frac{r^{N-1}}{C_{N-1}^{k-1}}p^{k}\left( r\right) h^{k}\left( w\left(
r\right) \right) ,\text{ }r=\left\vert x\right\vert .
\end{equation*}%
Then radial solutions of (\ref{6}) are any solution $w$ of the integral
equation%
\begin{equation*}
w\left( r\right) =1+\int_{0}^{r}\left( \frac{k}{t^{N-k}}\int_{0}^{t}\frac{%
s^{N-1}}{C_{N-1}^{k-1}}p^{k}\left( s\right) h^{k}\left( w\left( s\right)
\right) ds\right) ^{1/k}dt.\text{ }
\end{equation*}%
To establish a solution to this problem, we use successive approximation.
Define sequence $\left\{ w^{m}\right\} ^{m\geq 1}$ on $\left[ 0,\infty
\right) $ by%
\begin{equation*}
\left\{ 
\begin{array}{l}
w^{0}=1,\text{ }r\geq 0, \\ 
w^{m}\left( r\right) =1+\int_{0}^{r}\left( \frac{k}{t^{N-k}}\int_{0}^{t}%
\frac{s^{N-1}}{C_{N-1}^{k-1}}p^{k}\left( s\right) h^{k}\left( w^{m-1}\left(
s\right) \right) ds\right) ^{1/k}dt.%
\end{array}%
\right.
\end{equation*}%
We remark that, for all $r\geq 0$ and $m\in N$ 
\begin{equation*}
w^{m}\left( r\right) \geq 1\text{.}
\end{equation*}%
Moreover, proceeding by induction we conclude $\left\{ w^{m}\right\} ^{m\geq
1}$ are non-decreasing sequence on $\left[ 0,\infty \right) $. We note that $%
\left\{ w^{m}\right\} ^{m\geq 1}$ satisfies%
\begin{equation*}
\left\{ \frac{r^{N-k}}{k}\left[ \left( w^{^{m}}(r)\right) ^{\prime }\right]
^{k}\right\} ^{\prime }=\frac{r^{N-1}}{C_{N-1}^{k-1}}p^{k}\left( r\right)
h^{k}\left( w^{m-1}\left( r\right) \right) .
\end{equation*}%
By the monotonicity of $\left\{ w^{m}\right\} ^{m\geq 1}$ we have the
inequalities%
\begin{equation}
\left\{ \frac{r^{N-k}}{k}\left[ \left( w^{^{m}}(r)\right) ^{\prime }\right]
^{k}\right\} ^{\prime }=\frac{r^{N-1}}{C_{N-1}^{k-1}}p^{k}\left( r\right)
h^{k}\left( w^{m-1}\left( r\right) \right) \leq \frac{r^{N-1}}{C_{N-1}^{k-1}}%
p^{k}\left( r\right) h^{k}\left( w^{m}\left( r\right) \right) .  \label{8}
\end{equation}%
Choose $R>0$ so that $r^{N+\frac{N}{k}-2}p^{k}\left( r\right) $ are
non-decreasing for $r\geq R$. We are now ready to show that $w^{m}\left(
R\right) $ and $\left( w^{m}\left( R\right) \right) ^{\prime }$, both of
which are nonnegative, are bounded above independent of $m$. To do this, let 
\begin{equation*}
\phi ^{R}=\max \{p^{k}\left( r\right) :0\leq r\leq R\}\text{.}
\end{equation*}%
Using this and the fact that $\left( w^{m}\right) ^{\prime }\geq 0$, we note
that (\ref{8}) yields%
\begin{eqnarray*}
r^{N-k}\left[ \left( w^{^{m}}(r)\right) ^{\prime }\right] ^{k-1}\left(
w^{^{m}}(r)\right) ^{\prime \prime } &\leq &\frac{N-k}{k}r^{N-k-1}\left[
\left( w^{^{m}}(r)\right) ^{\prime }\right] ^{k}+r^{N-k}\left[ \left(
w^{^{m}}(r)\right) ^{\prime }\right] ^{k-1}\left( w^{^{m}}(r)\right)
^{\prime \prime } \\
&\leq &\phi _{1}^{R}\frac{r^{N-1}}{C_{N-1}^{k-1}}h^{k}\left( w^{m}\left(
r\right) \right) ,
\end{eqnarray*}%
and moreover%
\begin{equation*}
r^{N-k}\left[ \left( w^{^{m}}(r)\right) ^{\prime }\right] ^{k-1}\left(
w^{^{m}}(r)\right) ^{\prime \prime }\leq \phi ^{R}\frac{r^{N-k+k-1}}{%
C_{N-1}^{k-1}}h^{k}\left( w^{m}\left( r\right) \right) \leq R^{k-1}\phi ^{R}%
\frac{r^{N-k}}{C_{N-1}^{k-1}}h^{k}\left( w^{m}\left( r\right) \right) ,
\end{equation*}%
from which we have%
\begin{equation*}
\left[ \left( w^{^{m}}(r)\right) ^{\prime }\right] ^{k-1}\left(
w^{^{m}}(r)\right) ^{\prime \prime }\leq R^{k-1}\phi ^{R}\frac{1}{%
C_{N-1}^{k-1}}h^{k}\left( w^{m}\left( r\right) \right) .
\end{equation*}%
Multiply this by $\left( w^{^{m}}(r)\right) ^{\prime }$ we obtain%
\begin{equation}
\left\{ \left[ \left( w^{^{m}}(r)\right) ^{\prime }\right] ^{k+1}\right\}
^{\prime }\leq \frac{\left( k+1\right) R^{k-1}\phi ^{R}}{C_{N-1}^{k-1}}%
h^{k}\left( w^{m}\left( r\right) \right) \left( w^{^{m}}(r)\right) ^{\prime
}.  \label{ec}
\end{equation}%
Integrate (\ref{ec}) from $0$ to $r$ to get%
\begin{equation}
\left[ \left( w^{^{m}}(r)\right) ^{\prime }\right] ^{k+1}\leq \frac{\left(
k+1\right) R^{k-1}\phi ^{R}}{C_{N-1}^{k-1}}\int_{0}^{r}h^{k}\left(
w^{m}\left( s\right) \right) \left( w^{^{m}}(s)\right) ^{\prime }ds=\frac{%
\left( k+1\right) R^{k-1}\phi ^{R}}{C_{N-1}^{k-1}}%
\int_{1}^{w^{^{m}}(r)}h^{k}\left( s\right) ds  \label{9}
\end{equation}%
for $0\leq r\leq R$, which yields%
\begin{equation*}
\int_{1}^{w^{^{m}}(R)}\left[ \int_{1}^{t}h^{k}\left( s\right) ds\right]
^{-1/\left( k+1\right) }dt\leq \sqrt[k+1]{\frac{\left( k+1\right) \phi ^{R}}{%
C_{N-1}^{k-1}}}\cdot R^{\frac{2k}{k+1}}.
\end{equation*}%
It follows from the above relation and by the assumption (C2) that $%
w_{1}^{m}\left( R\right) $ is bounded above independent of $m$. Using this
fact in (\ref{9}) shows that the same is true of $\left( w^{m}\left(
R\right) \right) ^{\prime }$. Thus, the sequences $w^{m}\left( R\right) $
and $\left( w^{m}\left( R\right) \right) ^{\prime }$ are bounded above
independent of $m$.

Finally, we show that the non-decreasing sequences $w^{m}$ is bounded for
all $r\geq 0$ and all $m$. Multiplying the equation (\ref{8}) by $r^{N+\frac{%
N}{k}-2}\left( w^{^{m}}(r)\right) ^{\prime }$, we get%
\begin{equation}
\left\{ \left[ r^{\frac{N}{k}-1}\left( w^{m}\left( r\right) \right) ^{\prime
}\right] ^{k+1}\right\} ^{\prime }=\frac{k+1}{C_{N-1}^{k-1}}p^{k}\left(
r\right) h^{k}(w^{m}\left( r\right) )r^{N+\frac{N}{k}-2}\left( w^{m}\left(
r\right) \right) ^{\prime }.  \label{ma14}
\end{equation}%
Integrating from $R$ to $r$ gives%
\begin{equation*}
\left[ r^{\frac{N}{k}-1}\left( w^{m}\left( r\right) \right) ^{\prime }\right]
^{k+1}=\left[ R^{\frac{N}{k}-1}\left( w^{m}\left( R\right) \right) ^{\prime }%
\right] ^{k+1}+\frac{k+1}{C_{N-1}^{k-1}}\int_{R}^{r}p^{k}\left( s\right)
h^{k}(w^{m}\left( s\right) )s^{N+\frac{N}{k}-2}\left( w^{m}\left( s\right)
\right) ^{\prime }ds,
\end{equation*}%
for $r\geq R$. Noting that, by the monotonicity of $s^{N+\frac{N}{k}%
-2}p^{k}\left( s\right) $ for $r\geq s\geq R$, we get%
\begin{equation*}
\left[ r^{\frac{N}{k}-1}\left( w^{m}\left( r\right) \right) ^{\prime }\right]
^{k+1}\leq C+\frac{k+1}{C_{N-1}^{k-1}}r^{N+\frac{N}{k}-2}p^{k}\left(
r\right) H\left( w^{m}\left( r\right) \right) 
\end{equation*}%
where $C=\left[ R^{\frac{N}{k}-1}\left( w^{m}\left( R\right) \right)
^{\prime }\right] ^{k+1}$, which yields%
\begin{equation*}
r^{\frac{N}{k}-1}\left( w^{m}\left( r\right) \right) ^{\prime }\leq C^{\frac{%
1}{k+1}}+\left( \frac{k+1}{C_{N-1}^{k-1}}\right) ^{\frac{1}{k+1}}r^{\frac{N}{%
k}-\frac{2}{k+1}}p^{\frac{k}{k+1}}\left( r\right) H^{\frac{1}{k+1}}\left(
w^{m}\left( r\right) \right) 
\end{equation*}%
or, equivalently%
\begin{equation*}
\left( w^{m}\left( r\right) \right) ^{\prime }\leq C^{\frac{1}{k+1}}r^{1-%
\frac{N}{k}}+\left( \frac{k+1}{C_{N-1}^{k-1}}\right) ^{\frac{1}{k+1}}r^{1-%
\frac{2}{k+1}}p^{\frac{k}{k+1}}\left( r\right) H^{\frac{1}{k+1}}\left(
w^{m}\left( r\right) \right) 
\end{equation*}%
and hence%
\begin{equation}
\frac{d}{dr}\int_{w^{m}\left( R\right) }^{w^{m}\left( r\right) }\left[
H\left( t\right) \right] ^{-1/\left( k+1\right) }dt\leq C^{\frac{1}{k+1}%
}r^{1-\frac{N}{k}}H^{-\frac{1}{k+1}}\left( w^{m}\left( r\right) \right)
+\left( \frac{k+1}{C_{N-1}^{k-1}}\right) ^{\frac{1}{k+1}}\left(
r^{k-1}p^{k}\left( r\right) \right) ^{\frac{1}{k+1}}.  \label{in}
\end{equation}%
Inequality (\ref{in}) combined with%
\begin{eqnarray*}
\frac{1}{\sqrt{2}}\sqrt{2\cdot \left( s^{k-1}p^{k}\left( s\right) \right) ^{%
\frac{2}{k+1}}} &=&\frac{1}{\sqrt{2}}\sqrt{2\cdot s^{\frac{1+\varepsilon }{2}%
}\left( s^{k-1}p^{k}\left( s\right) \right) ^{\frac{2}{k+1}}s^{\frac{%
1-\varepsilon }{2}}} \\
&\leq &\frac{1}{\sqrt{2}}\left[ s^{1+\varepsilon }\left( s^{k-1}p^{k}\left(
s\right) \right) ^{\frac{2}{k+1}}+s^{-1-\varepsilon }\right] 
\end{eqnarray*}%
gives%
\begin{equation*}
\begin{array}{ll}
\int_{w^{m}\left( R\right) }^{w^{m}\left( r\right) }\left[ H\left( t\right) %
\right] ^{-1/\left( k+1\right) }dt & \leq C^{\frac{1}{k+1}}\int_{R}^{r}t^{1-%
\frac{N}{k}}H^{-\frac{1}{k+1}}\left( w^{m}\left( t\right) \right) dt \\ 
& +\frac{1}{\sqrt{2}}\left( \frac{k+1}{C_{N-1}^{k-1}}\right) ^{\frac{1}{k+1}}%
\left[ \int_{R}^{r}t^{1+\varepsilon +\frac{2\left( k-1\right) }{k+1}}\left(
p\left( t\right) \right) ^{\frac{2k}{k+1}}dt+\int_{R}^{r}t^{-1-\varepsilon
}dt\right]  \\ 
& \leq C^{\frac{1}{k+1}}H^{-\frac{1}{k+1}}\left( w^{m}\left( R\right)
\right) \int_{R}^{r}t^{1-\frac{N}{k}}dt \\ 
& +\frac{1}{\sqrt{2}}\left( \frac{k+1}{C_{N-1}^{k-1}}\right) ^{\frac{1}{k+1}}%
\left[ \int_{R}^{r}t^{1+\varepsilon +\frac{2\left( k-1\right) }{k+1}}\left(
p\left( t\right) \right) ^{\frac{2k}{k+1}}dt+\frac{1}{\varepsilon
R^{\varepsilon }}\right] \text{.}%
\end{array}%
\end{equation*}%
The above relation is needed in proving the bounded of the function $\left\{
w^{m}\right\} ^{m\geq 1}$ in the following. Indeed, since for each $%
\varepsilon >0$ the right side of this inequality is bounded independent of $%
m$ (note that $w^{m}\left( t\right) \geq 1$), so is the left side and hence,
in light of (C2), the sequence $\left\{ w^{m}\right\} ^{m\geq 1}$ is a
bounded sequence and so $\left\{ w^{m}\right\} ^{m\geq 1}$ are bounded
sequence. Thus $\left\{ w^{m}\right\} ^{m\geq 1}\rightarrow w$ as $%
m\rightarrow \infty $ and the limit functions $w$ are positive entire
bounded solutions of equation (\ref{6}).

\subparagraph{Proof of the Theorem \protect\ref{th2}.}

We know that for any $a_{1}>0$ a solution of%
\begin{equation*}
v\left( r\right) =a_{1}+\int_{0}^{r}\left( \frac{k}{t^{N-k}}\int_{0}^{t}%
\frac{s^{N-1}}{C_{N-1}^{k-1}}p^{k}\left( s\right) h^{k}\left( v\left(
s\right) \right) ds\right) ^{1/k}dt,
\end{equation*}%
exists, at least, small $r$. Since $v^{\prime }\geq 0$, the only way that
the solution can become singular at $R$ is for $v\left( r\right) \rightarrow
\infty $ as $r\longrightarrow \infty $. Thus, we can show that, for each $R>0
$, there exists $C_{R}>0$ so that $v\left( R\right) \leq C_{R}$, we have
existence. To this end, let $M_{R}=\max \left\{ p\left( r\right) \left\vert
0\leq r\leq R\right. \right\} $ and consider the equation%
\begin{equation*}
w\left( r\right) =a_{2}+M_{R}\int_{0}^{r}\left( \frac{k}{t^{N-k}}\int_{0}^{t}%
\frac{s^{N-1}}{C_{N-1}^{k-1}}h^{k}\left( v\left( s\right) \right) ds\right)
^{1/k}dt
\end{equation*}%
where $a_{2}>a_{1}$. We next observe that the solution to this equation
exists for all $r\geq 0$ and of course, it is a solution to $%
S_{k}^{1/k}\left( \lambda \left( D^{2}w\left( r\right) \right) \right)
=M_{R}h\left( w\right) $ on $\mathbb{R}^{N}$ which is treated in \cite[%
(Theorem 1.1, p. 177)]{BAO}. We now show that $v\left( r\right) \leq w\left(
r\right) $ for all $0\leq r\leq R$ and hence we conclude the proof of
existence. Clearly $v\left( 0\right) <w\left( 0\right) $ so that $v\left(
r\right) <w\left( r\right) $ for at least all $r$ near zero. Let 
\begin{equation*}
r_{0}=\sup \left\{ r\left\vert v\left( s\right) <w\left( s\right) \text{ for
all }s\in \left[ 0,r\right] \right. \right\} \text{.}
\end{equation*}%
If $r_{0}=R$, then we are done. Thus assume that $r_{0}<R$. It follows from
assumption $a_{2}>a_{1}$ that 
\begin{eqnarray*}
v\left( r_{0}\right)  &=&a_{1}+\int_{0}^{r_{0}}\left( \frac{k}{t^{N-k}}%
\int_{0}^{t}\frac{s^{N-1}}{C_{N-1}^{k-1}}p^{k}\left( s\right) h^{k}\left(
v\left( s\right) \right) ds\right) ^{1/k}dt \\
&<&a_{2}+M_{R}\int_{0}^{r_{0}}\left( \frac{k}{t^{N-k}}\int_{0}^{t}\frac{%
s^{N-1}}{C_{N-1}^{k-1}}h^{k}\left( v\left( s\right) \right) ds\right)
^{1/k}dt=w\left( r_{0}\right) .
\end{eqnarray*}%
Thus there exists $\varepsilon >0$ so that $v\left( r\right) <w\left(
r\right) $ for all $\left[ 0,r+\varepsilon \right) $, contradicting the
definition of $r_{0}$. Thus we conclude that $v<w$ on $\left[ 0,R\right] $
for all $R>0$ and hence $v$ is a nontrivial entire solution of (\ref{11}).
Now let $u$ be any nonnegative nontrivial entire solution of (\ref{11}) and
suppose $p$ satisfies 
\begin{equation*}
\int_{0}^{\infty }\left( \frac{k}{t^{N-k}}\int_{0}^{t}\frac{s^{N-1}}{%
C_{N-1}^{k-1}}p^{k}\left( s\right) ds\right) ^{1/k}dt=\infty .
\end{equation*}%
Since $u$ is nontrivial and non-negative, there exists $R>0$ so that $%
u\left( R\right) >0$. On the other hand since $u^{\prime }\geq 0$, we get $%
u\left( r\right) \geq u\left( R\right) $ for $r\geq R$ and thus from%
\begin{equation*}
u\left( r\right) =u\left( 0\right) +\int_{0}^{r}\left( \frac{k}{t^{N-k}}%
\int_{0}^{t}\frac{s^{N-1}}{C_{N-1}^{k-1}}p^{k}\left( s\right) h^{k}\left(
u\left( s\right) \right) ds\right) ^{1/k}dt,\text{ }
\end{equation*}%
since $u$ will satisfy that equation for all $r\geq 0,$ we get 
\begin{equation*}
\begin{array}{l}
u\left( r\right) =u\left( 0\right) +\int_{0}^{r}\left( \frac{k}{t^{N-k}}%
\int_{0}^{t}\frac{s^{N-1}}{C_{N-1}^{k-1}}p^{k}\left( s\right) h^{k}\left(
u\left( s\right) \right) ds\right) ^{1/k}dt \\ 
\text{ \ \ \ \ \ \ }\geq u\left( R\right) +h\left( u\left( R\right) \right)
\int_{R}^{r}\left( \frac{k}{t^{N-k}}\int_{R}^{t}\frac{s^{N-1}}{C_{N-1}^{k-1}}%
p^{k}\left( s\right) ds\right) ^{1/k}dt\rightarrow \infty \text{ as }%
r\rightarrow \infty .%
\end{array}%
\end{equation*}%
Conversely, assume that\ $h$ satisfy (C1), (C3) and that $w$ is a
nonnegative entire large solution of (\ref{11}). Note also, that $w$
satisfies%
\begin{equation*}
\left[ \frac{r^{N-k}}{k}\left( w^{^{\prime }}(r)\right) ^{k}\right] ^{\prime
}=\frac{r^{N-1}}{C_{N-1}^{k-1}}p^{k}\left( r\right) h^{k}\left( w\left(
r\right) \right) .
\end{equation*}%
Using the monotonicity of \textit{\ }$r^{N+\frac{N}{k}-2}p\left( r\right) $%
\textit{\ } we can apply similar arguments used in obtaining Theorem \ref%
{th1} to get%
\begin{equation*}
\left( w\left( r\right) \right) ^{\prime }\leq C^{\frac{1}{k+1}}r^{1-\frac{N%
}{k}}+\left( \frac{k+1}{C_{N-1}^{k-1}}\right) ^{\frac{1}{k+1}}r^{1-\frac{2}{%
k+1}}p^{\frac{k}{k+1}}\left( r\right) H^{\frac{1}{k+1}}\left( w\left(
r\right) \right) ,
\end{equation*}%
which we may rewrite as%
\begin{equation}
\begin{array}{ll}
\int_{w\left( R\right) }^{w\left( r\right) }\left[ H\left( t\right) \right]
^{-1/\left( k+1\right) }dt & \leq C^{\frac{1}{k+1}}\int_{R}^{r}t^{1-\frac{N}{%
k}}H^{-\frac{1}{k+1}}\left( w\left( t\right) \right) dt \\ 
& +\frac{1}{\sqrt{2}}\left[ \int_{R}^{r}t^{1+\varepsilon +\frac{2\left(
k-1\right) }{k+1}}\left( p\left( t\right) \right) ^{\frac{2k}{k+1}%
}dt+\int_{R}^{r}t^{-1-\varepsilon }dt\right]  \\ 
& \leq C^{\frac{1}{k+1}}H^{-\frac{1}{k+1}}\left( w\left( R\right) \right)
\int_{R}^{r}t^{1-\frac{N}{k}}dt \\ 
& +\frac{1}{\sqrt{2}}\left( \frac{k+1}{C_{N-1}^{k-1}}\right) ^{\frac{1}{k+1}}%
\left[ \int_{R}^{r}t^{1+\varepsilon +\frac{2\left( k-1\right) }{k+1}}\left(
p\left( t\right) \right) ^{\frac{2k}{k+1}}dt+\frac{1}{\varepsilon
R^{\varepsilon }}\right]  \\ 
& \leq C_{R}\int_{R}^{r}t^{1-\frac{N}{k}}dt+\frac{1}{\sqrt{2}}\left( \frac{%
k+1}{C_{N-1}^{k-1}}\right) ^{\frac{1}{k+1}}\int_{R}^{r}t^{1+\varepsilon +%
\frac{2\left( k-1\right) }{k+1}}\left( p\left( t\right) \right) ^{\frac{2k}{%
k+1}}dt%
\end{array}
\label{final}
\end{equation}%
where 
\begin{equation*}
C_{R}=C^{\frac{1}{k+1}}H^{-\frac{1}{k+1}}\left( w\left( R\right) \right) +%
\frac{1}{\sqrt{2}}\frac{1}{\varepsilon R^{\varepsilon }}\left( \frac{k+1}{%
C_{N-1}^{k-1}}\right) ^{\frac{1}{k+1}}.
\end{equation*}%
By taking $r\rightarrow \infty $ in (\ref{final}) we obtain (\ref{13}) since 
$w$ is large and $h$ satisfies (C3). These observations completes the proof
of the theorem.

\subparagraph{\textbf{Proof of the Theorem \protect\ref{th3} and \protect\ref%
{th4}.}}

In order, to obtain the conclusion, combine the proof of \textbf{Theorem \ref%
{th1} }and\textbf{\ \ref{th2}} with some technical results from \cite{COV}
and \cite{COVS}.

\qed

\end{document}